\PassOptionsToPackage{dvipsnames}{xcolor}
\documentclass[final]{siamart251216}

\usepackage{graphicx} % Required for inserting images
\usepackage{todonotes} % For orange notes
\usepackage[utf8]{inputenc} % Fonts
\usepackage{algorithmic}

\usepackage{mathtools}
\usepackage{stmaryrd} 

\usepackage{amsfonts}
\usepackage{amssymb}

\newsiamremark{remark}{Remark}
\newsiamremark{example}{Example}

\DeclareMathOperator{\trace}{Tr}

\DeclareMathOperator{\diag}{diag}
\DeclareMathOperator{\orth}{orth}	% orthogonal basis
\DeclareMathOperator{\vspan}{span} 	% vector span
\DeclarePairedDelimiter{\norm}{\lVert}{\rVert}
\DeclarePairedDelimiter{\abs}{\lvert}{\rvert}

\newcommand{\R}{\mathbb{R}} % Real numbers
\newcommand{\C}{\mathbb{C}} % Complex numbers
\renewcommand{\vec}[1]{\boldsymbol{#1}}		% for bold vectors
\newcommand{\de}{\mathrm{d}} % The "d" in integrals
\newcommand{\K}{\mathcal{K}} % Generic Krylov subspace
\newcommand{\kryl}{\K}

\newcommand{\poly}{\mathbb{P}}
\renewcommand{\b}{\vec{b}}
\newcommand{\B}{\vec{B}}
\newcommand{\V}{\vec{V}}
\let\textcaron\v % Preserve the text accent used in author names.
\renewcommand{\v}{\ifmmode\vec{v}\else\expandafter\textcaron\fi}
\newcommand{\W}{\vec{W}}

\renewcommand{\u}{\vec{u}}
\newcommand{\X}{\vec{X}}
\newcommand{\x}{\vec{x}}
\newcommand{\Y}{\vec{Y}}

\newcommand{\E}{\vec{E}}
\newcommand{\e}{\vec{e}}

\title{Global Lanczos with stochastic error correction for trace estimation of matrix functions}

\author{Michele Rinelli\thanks{KU Leuven, Leuven, Belgium
  (\email{michele.rinelli@kuleuven.be}, \email{robbe.vermeiren@kuleuven.be}).}
  \and Igor Simunec\thanks{Aalto University, Espoo, Finland (\email{igor.simunec@aalto.fi}).}
  \and Robbe Vermeiren\footnotemark[1]}

\headers{Trace estimation with global Lanczos}{M. Rinelli, I. Simunec, and R. Vermeiren}

\usepackage[caption=false]{subfig} % subfigures, as recommended for the SIAM class

\numberwithin{equation}{section}

\begin{document}

\maketitle

\begin{abstract}
	We propose a new method for computing the trace of a matrix function $f(A)$, where $A$ is a large, sparse, symmetric matrix. 
	Our approach combines a deterministic global Lanczos approximation with a stochastic correction of its error.
Starting the global Lanczos process with the identity matrix yields a direct approximation of the trace through a small tridiagonal matrix. The resulting basis matrices are polynomials in $A$, allowing the process to be implemented using sparse matrix-matrix multiplications when low-degree polynomials in $A$ have moderate fill-in. However, increased fill-in may limit the number of iterations that can be performed efficiently. 
To address this limitation, we derive an exact representation of the global Lanczos error as the trace of a matrix function and apply a stochastic estimator to the remainder. Even a moderately accurate estimate can substantially improve the accuracy, combining the deterministic method with the rapid initial convergence of stochastic trace estimators. We discuss error bounds and variance reduction, and present numerical experiments that compare the resulting method with existing ones.
\end{abstract}

\begin{keywords}
trace estimation, matrix functions, global Lanczos, stochastic trace estimation
\end{keywords}

\begin{MSCcodes}
65F60, 65F50, 65C05
\end{MSCcodes}

\section{Introduction}

We consider the approximation of the trace of a matrix function $f(A)$, where $A\in\mathbb{R}^{n\times n}$ is a large, sparse, symmetric matrix and $f$ is a scalar function defined on its spectrum. Denoting the eigenvalues of $A$ by $\lambda_1,\ldots,\lambda_n$, the quantity of interest is
\begin{equation*}
\trace(f(A))=\sum_{i=1}^{n}f(\lambda_i).
\end{equation*}
This problem has several applications, including the log-determinant computation~\cite{CortinovisKressner22}, the von Neumann entropy in quantum statistical mechanics~\cite{LandauLifshitz1958,Wehrl1978} and total $\pi$-electron energies of alternant hydrocarbons in H\"uckel theory~\cite{CoulsonOLearyMallion1978}. In statistics, it enters Gaussian likelihood evaluation~\cite{aune2014parameter}. The trace of an inverse matrix, i.e., $f(x)=x^{-1}$, arises in lattice quantum chromodynamics, generally for nonsymmetric matrices~\cite{StathopoulosLaeuchliOrginos2013}. In graph theory, traces of functions of adjacency matrices or graph Laplacians define quantities such as graph energy~\cite{EstradaBenzi2017}, von Neumann graph entropy~\cite{BraunsteinGhoshSeverini2006}, and centrality measures like the Estrada index and total resolvent subgraph centrality~\cite{EstradaHigham2010}. For large matrices, computing the full spectrum or explicitly forming $f(A)$ can be prohibitively expensive, so it is often much cheaper to rely on approximations.

Stochastic trace estimators approximate the trace of a matrix $M$ using matrix-vector products or quadratic forms. The classical Monte Carlo estimator of Girard~\cite{Girard1989} and Hutchinson~\cite{Hutchinson89}, commonly known as Hutchinson's estimator,
 averages quadratic forms
 $\vec x^T M \vec x$
 over independent random vectors satisfying $\mathbb{E}[\vec x\vec x^T]=I$. 
\space
Randomized low-rank approaches to trace estimation were studied by Saibaba, Alexanderian, and Ipsen~\cite{SaibabaLowrank17}. Meyer, Musco, Musco, and Woodruff~\cite{MeyerMuscoMuscoWoodruff2021} developed the \texttt{Hutch++} method, which computes the trace of a randomized low-rank approximation and applies Hutchinson's estimator to the remainder. \space Inspired by \texttt{Hutch++}, Epperly, Tropp, and Webber~\cite{EpperlyTroppWebber24} designed \texttt{XTrace} and \texttt{XNysTrace} which also use randomized low-rank approximations but exploits \emph{exchangeability} to reuse the random samples more efficiently. When $M=f(A)$, the matrix-vector products and quadratic forms are usually approximated using Krylov methods; see, for example, Saad~\cite{Saad1992} and G\"uttel~\cite{GuttelRatKryl13} for polynomial and rational Krylov, respectively, and Benzi, Rinelli, and Simunec~\cite{BenziRinelliSimunec2023} for their use in trace estimation. Chen and Hallman~\cite{ChenHallman23} developed a Krylov-aware estimator that integrates the low-rank approximation and stochastic estimation with the underlying Krylov computations.

Probing methods, developed by Tang and Saad~\cite{tang2012probing} and analyzed for trace estimation by Frommer, Schimmel, and Schweitzer~\cite{FrommerSchimmelSchweitzer21}, provide another approach to approximate $\operatorname{tr}(f(A))$ for a sparse matrix $A$. They use a distance-$d$ coloring of the graph associated with $A$ to construct probing vectors and approximate the trace by summing the corresponding quadratic forms. Probing is effective when the entries of $f(A)$ decay rapidly with graph distance. A stochastic variant, in which the nonzero entries of the probing vectors are replaced by random signs, was proposed by Aune, Simpson, and Eidsvik~\cite{aune2014parameter}, while Frommer, Rinelli, and Schweitzer~\cite{FrommerRinelliSchweitzer25} provided a convergence analysis.

A related problem is the computation of 
\begin{equation*}
\trace(\B^T f(A)\B)=\sum_{j=1}^s \b_j^T f(A)\b_j,
\end{equation*}
where $\vec B=[\b_1,\dots,\b_s]$ is a block vector. 
The global Lanczos algorithm approximates the entire sum through a single process and a small tridiagonal matrix, without needing to run separate single-vector projection methods. This approach was studied, e.g., by Bellalij, Reichel, Rodriguez, and Sadok~\cite{ReichelGlobal15}. The connection with Gauss and Gauss--Radau quadrature, global Lanczos also yields upper and lower bounds under suitable sign conditions on the derivatives of $f$.

In this paper, we use the global Lanczos algorithm with $\B=I_n$ to approximate $\trace(f(A))$ directly, and we refer to this choice as the \emph{full-size approach}. \space This choice was discussed by D\'iaz Fuentes et al.~\cite{FDFM23} and used in a matrix-product-operator representation by August, Ba\~nuls, and Huckle~\cite{AugustBanulsHuckle2017}. Here, we assume that $A$ is explicitly available, allowing us to exploit the sparsity of the full-size basis matrices, which are polynomials in $A$. Provided that low-degree polynomials in $A$ have moderate fill-in, the iterations can be performed efficiently using sparse matrix--matrix multiplications. As fill-in increases, however, further iterations may become too expensive before the approximation reaches the desired accuracy.

To improve this approximation, we derive an exact error representation based on Hermite interpolation. Let $\psi_m^{GL}$ be the approximation after $m$ Lanczos steps and $h_m$ the Hermite interpolant of $f$ at the Lanczos nodes. Let $M = f(A) - h_m(A)$, then we show that
\begin{equation*}
\trace(f(A)) - \psi_m^{GL} = \trace(M).
\end{equation*}
This identity motivates our \emph{stochastic error correction}: we add a stochastic estimate of $\trace(M)$, obtained with Hutchinson, XTrace, or XNysTrace, to the deterministic approximation. Our analysis provides bounds on the Frobenius norm of the error matrix $r_m(A)$, which controls the sampling variance. A smaller error matrix can therefore allow an accurate corrected approximation with fewer random samples. In numerical experiments, we also observe examples where XTrace/XNysTrace are more effective on the error matrix than on $f(A)$, even when the latter does not admit an accurate low-rank approximation.

The remainder of the paper is organized as follows. In~\cref{sec:global-Lanczos}, we review the global Lanczos process, discuss the full-size approach and its applicability, present error bounds and derive an exact error representation. In~\cref{sec:other-trace-estimators}, we review stochastic and probing estimators. In~\cref{sec:stochastic-error-correction}, we develop the stochastic error correction, analyze its theoretical properties, and discuss implementation aspects. Finally, in~\cref{sec:numerical-experiments}, we compare the proposed method with existing trace estimators through numerical experiments.

\subsection{Notation}
We denote matrices by $A \in \R^{n \times n}$, vectors by $\vec b \in \R^{n}$, and block vectors by $\vec B \in \R^{n \times \ell}$. For symmetric $A$, $\sigma(A)$ denotes its spectrum, and $\lambda_{\min}$ and $\lambda_{\max}$ its smallest and largest eigenvalues. We denote by $\norm{\cdot}_F$ the Frobenius norm and by $\norm{\cdot}_2$ the Euclidean norm of a vector or its induced matrix norm. The Frobenius inner product is $\langle\V,\W\rangle=\trace(\V^T\W)$; in the full-size case, both arguments are $n\times n$ matrices. For square $M$, $\diag(M)$ is the diagonal matrix with diagonal entries $M_{ii}$, and $M_{\mathrm{off}}=M-\diag(M)$ denotes its off-diagonal part.

We write $f[x_0,\ldots,x_k]$ for divided differences of $f$~\cite{deBoor05}, with repeated nodes interpreted through derivatives. The symbols $\mathbb{E}$, $\mathbb{V}$, and $\mathbb{P}$ denote expectation, variance, and probability, respectively. We also use the shorthand $\mathcal{I}[f,\B]:=\trace(\B^T f(A)\B)$.

\section{Trace estimation with global block Lanczos}
\label{sec:global-Lanczos}
\subsection{Krylov subspaces and Lanczos}
Given a matrix $A\in\R^{n\times n}$ and a vector $\b\in\R^{n}$, the Krylov subspace of order $m$ associated with the pair $(A,\b)$ is 
\begin{equation*}
	\kryl_m(A,\b) = \vspan(\b,A\b,\dots, A^{m-1}\b).
\end{equation*}
Krylov subspaces see plenty of usage in many problems in numerical linear algebra and applied sciences, such as linear system solving, matrix function computations, and model order reduction; see, e.g.,~\cite{Higham08,LiesenStrakosKrylov13}.

\subsection{Global block Lanczos}
The $m$-th order global block Krylov subspace \cite{FLS17} associated with $A\in\R^{n\times n}$ and a block vector $\B\in\R^{n\times \ell}$ is defined as
\begin{equation*}
    \kryl_m(A,\B) := \left\{ \sum_{k=0}^{m-1} c_kA^k\B: c_k\in\R \right\}.
\end{equation*}
Define the following inner product on block vectors
\begin{equation}
	\langle\V,\W\rangle := \trace(\V^T\W).
\end{equation}
Similar to the classical Lanczos algorithm~\cite{Saad2003}, $m$ steps of the global block Lanczos algorithm compute a nested basis $\vec V_1, \dots, \vec V_m$ of $\kryl_m(A,\B)$, with $\vec V_k \in \R^{n \times \ell}$, which is orthonormal with respect to $\langle \cdot, \cdot \rangle$, i.e., $\trace(\vec V_i^T\vec V_j) = \delta_{i,j}$.
Note that $\langle \X,\X\rangle = \|\X\|_F^2$ for $\X\in\R^{n\times \ell}$, so $\X$ can be normalized by taking $\Y = \X/\beta$ with $\beta = \norm{\X}_F$, so that $\langle\Y,\Y\rangle = 1$.

\begin{algorithm}[t]
\caption{Global Lanczos}
\label{alg:global-lanczos}
\begin{algorithmic}[1]
\REQUIRE Symmetric $A\in\R^{n\times n}$, $\B \in\R^{n\times \ell}$, function $f$, number of steps $m$
\ENSURE $\psi_m^{GL}[f,\B] \approx \mathcal{I}[f,\B]$, $T_m$

\STATE \textbf{Initialize.} $\V_0 = \vec 0$, $\beta_1 = \|\B\|_F$, $\V_1 = \B / \beta_1$.

\FOR{$k=1,2,\dots,m$}
    \STATE $\vec X = A\V_k-\beta_k\V_{k-1}$, $\alpha_k = \langle\V_k,\X\rangle$
    \STATE $\X = \X-\alpha_k\V_k$
    \STATE $\beta_{k+1} = \|\X\|_F$, $\V_{k+1} = \X / \beta_{k+1}$
\ENDFOR
\STATE $\psi_m^{GL}[f,\B] = \beta_1^2\, \e_1^T f(T_m)\, \e_1$
\end{algorithmic}
\end{algorithm}

The orthonormal basis vectors satisfy a three-term recurrence relation:
\begin{equation}
	\label{eqn:global-lanczos-recurrence}
    A\V_k = \beta_{k+1}\V_{k+1} + \alpha_k\V_k + \beta_k\V_{k-1},\qquad k = 1,\dots,m.
\end{equation}
The main distinction from classical block Lanczos is that the recurrence coefficients are scalars, and this imposes additional structure on the recurrence matrix.
In order to write this in matrix form, define the $m \times m$ tridiagonal matrix
\begin{equation*}
    T_m = \begin{bmatrix}
\alpha_1 & \beta_2\\
\beta_2 & \alpha_2 & \ddots\\
& \ddots & \ddots & \beta_{m}\\
& & \beta_m & \alpha_m
\end{bmatrix},
\end{equation*}
and the block tridiagonal matrix $\mathcal{T}_m = T_m \otimes I_\ell \in\R^{\ell m\times \ell m}$. Finally, set
    $\mathcal{V}_m = [\V_1,\dots,\V_m],\quad \vec E_m = \vec e_m\otimes I_\ell \in\C^{m\ell\times \ell},$ 
with $\e_m$ denoting the $m$-th vector of the canonical basis in $\R^m$. Then,
\begin{equation*}
    A\mathcal{V}_m = \mathcal{V}_m \mathcal{T}_m + \beta_{m+1}\V_{m+1}\vec \E_m^T.
\end{equation*}
It is worth observing that the block vectors $\V_k$ are only orthonormal with respect to $\langle\cdot,\cdot\rangle$, so $\mathcal{V}_m^T\mathcal{V}_m\neq I_{m\ell}$ and $\mathcal{V}_m^T A \mathcal{V}_m \neq \mathcal{T}_m$.

\subsection{Gauss quadrature for $\mathcal{I}[f,\B]$} 
Consider the spectral decomposition $A = \sum_{i=1}^n \lambda_i \vec z_i \vec z_i^T$, where $(\lambda_i, \vec z_i)$ are the eigenpairs of $A$, then
\begin{equation}
    \label{eq:trace-functional-integral-representation}
\mathcal{I}[f,\B] = \sum_{i=1}^n f(\lambda_i) \mu_i =\int_{\R} f(x) \,\de\mu(x),
\end{equation}
where $\mu$ is the discrete measure\footnote{$\delta_{\lambda_i}$ is the Dirac delta measure concentrated in $\lambda_i$, i.e.\ $\delta_{\lambda_i}(\Omega) = 1$ if $\lambda_i\in\Omega$ and $0$ otherwise.} defined on $\Omega\subset \R$ by $\mu(\Omega) = \sum_{i=1}^n \delta_{\lambda_i}(\Omega)\,\mu_i$, with positive \textit{weights} $\mu_i = \trace(\B^T\vec z_i\vec z_i^T\B) = \vec z_i^T\B\B^T\vec z_i$, \textit{nodes} $\lambda_i$ and total mass $\norm{\vec B}_F^2$.
Since $\mu$ is supported on $\sigma(A)$, the integral can be taken over any $[a,b]\supseteq\sigma(A)$.

In view of~\cref{eq:trace-functional-integral-representation}, approximating $\mathcal{I}[f,\B]$ amounts to approximating an integral, so it is natural to apply a quadrature rule. Bellalij et al.~\cite{ReichelGlobal15} use the \emph{global Lanczos approximation}
\begin{equation}\label{eq:GL-approximation}
\psi_m^{GL}[f, \vec B] := \beta_1^2\, \e_1^T f(T_m)\, \e_1 = \|\B\|_F^2\, \e_1^T f(T_m)\, \e_1,
\end{equation}
which is computed by~\cref{alg:global-lanczos}, and show that this is precisely the optimal Gauss quadrature rule for the integral~\eqref{eq:trace-functional-integral-representation}. This intimate connection between Krylov subspaces and Gauss quadrature is well known; see, for example, the book by Golub and Meurant~\cite{GolubMeurantBook10}.

\begin{theorem}[{\cite{ReichelGlobal15}}]\label{thm:gauss-approximation}
    Consider the spectral decomposition $T_m=\sum_{k=1}^m \theta_k \u_k\u_k^T$. Then, 
\begin{equation}\label{eq:gauss-quadrature-rule}
\psi^{GL}_m[f, \vec B] = \beta_1^2 \sum_{k=1}^m w_kf(\theta_k),\quad w_k = |\e_1^T\vec u_k|^2.
\end{equation}
    Moreover, if $p_0,\dots,p_m$ are the polynomials of degree $\deg(p_k)=k$ defined by
	$\V_k = p_{k-1}(A)\B$, with $\V_k$ from~\cref{alg:global-lanczos}, then they are
	orthonormal with respect to $\mu$, i.e.\ $\int_{\mathbb{R}} p_i p_j \,\de\mu
	= \delta_{i,j}$, and the roots of $p_m$ are $\theta_1,\dots,\theta_m$. Finally,
	$\psi^{GL}_m[f, \vec B]$ is the $m$-th Gauss quadrature rule for $\mu$, so it is exact for any polynomial $f$ of degree at most $2m-1$.
\end{theorem}

A direct consequence, which to our knowledge has not been stated explicitly for the global Lanczos approximation, is a bound based on polynomial approximations of $f$ that reveals the rate of convergence. A similar result for classical block Lanczos is given by \cite[Lemma~2.2]{ChenHallman23}.
\begin{corollary}
    \label{cor:convergence-Gauss-quadrature}
    Let $a<b$ be such that $\sigma(A)\subset[a,b]$, and let $f$ be continuous over $[a,b]$. Then
    \begin{equation*}
        |\mathcal{I}[f,\B] - \psi^{GL}_m[f, \vec B]| \leq 2\|\B\|_F^2 \, E_{2m-1}(f, [a, b]),
    \end{equation*}
	where 
	\begin{equation*}
		E_{2m-1}(f, [a, b]) := \min_{p \in \poly_{2m-1}} \max_{x\in[a,b]} \abs{f(x)-p(x)}.
	\end{equation*}
\end{corollary}
\begin{proof}
    Let $p\in\poly_{2m-1}$. By polynomial exactness, $\mathcal{I}[p,\B] = \psi^{GL}_m[p, \vec B]$, so by linearity
    \begin{equation*}
        |\mathcal{I}[f,\B] - \psi^{GL}_m[f, \vec B]| \leq |\mathcal{I}[f-p,\B]| + \|\B\|_F^2\, |\e_1^T(f-p)(T_m)\e_1|.
    \end{equation*}
    Writing $\B = [\b_1,\dots,\b_\ell]$, the first term satisfies
    \begin{equation*}
        |\mathcal{I}[f-p,\B]| \leq \sum_{j=1}^\ell |\b_j^T(f-p)(A)\b_j| \leq \|\B\|_F^2\|(f-p)(A)\|_2 = \|\B\|_F^2 \max_{x\in\sigma(A)}|f(x)-p(x)|,
    \end{equation*}
    and the second is at most $\|\B\|_F^2\|(f-p)(T_m)\|_2 = \|\B\|_F^2\max_{x\in\sigma(T_m)}|f(x)-p(x)|$. Since $\sigma(A),\sigma(T_m)\subset[a,b]$, both are bounded by $\|\B\|_F^2\max_{x\in[a,b]}|f(x)-p(x)|$, and minimizing over $p$ concludes the proof.
\end{proof}

\subsection{Full-size global block Lanczos for trace estimation}
\label{subsec:full-size-global-lanczos}

\cref{alg:global-lanczos} can be used to approximate $\trace(f(A))$ as follows \cite{ReichelGlobal15}. Let $I_n = [\vec B_1, \dots, \vec B_s]$ be a partition of the columns of $I_n$, where we assume for simplicity that $\vec B_i \in \R^{n \times \ell}$ for all $i$. Then, we have the approximation 
\begin{equation}
	\label{eqn:global-lanczos-approx-trace}
	\trace(f(A)) = \sum_{i = 1}^s \mathcal{I}[f, \vec B_i] \approx \sum_{i = 1}^s \psi_m^{GL}[f, \vec B_i]. 
\end{equation}
In \cite{ReichelGlobal15}, only relatively small block sizes were considered, with $\ell \le 100$.

We now set $s = 1$ and $\B = I_n$ in this construction. Then
\begin{equation} \label{eq:full-size-GL}
    \mathcal{I}[f,I_n] = \trace(f(A)) = n\,\e_1^T f(T_m)\e_1,
\end{equation}
and~\cref{eq:GL-approximation} and \cref{eqn:global-lanczos-approx-trace} become $\trace(f(A)) \approx \psi_m^{GL}[f, I_n] = n\,\e_1^T f(T_m)\e_1$.
The initialization in~\cref{alg:global-lanczos} is $\beta_1 = \|I_n\|_F = \sqrt n$ and $\V_1 = I_n/\sqrt n$. It follows from \cref{thm:gauss-approximation} that $\vec V_j = p_{j-1}(A)$ for $j = 1, \dots, m+1$. In the following, we will use the simplified notation $\psi_m^{GL} := \psi_m^{GL}[f, I_n]$ when there is no ambiguity.

The main computational constraint is the storage and manipulation of the matrices $\V_j$. Although each step uses the same three-term recurrence as~\cref{alg:global-lanczos}, the main advantage of the full-size formulation is the possibility of using sparse matrix-matrix operations. Indeed, the products $A\V_j$ and the Frobenius inner products are matrix-matrix products of $n \times n$ matrices, which can be computed efficiently if all the matrices $p_{j-1}(A)$ are sufficiently sparse. The efficiency of the method therefore depends on the sparsity of $A$ and on the fill-in of low-degree polynomials in $A$.

This fill-in can be expressed in terms of the undirected graph $G(A)$ associated with the matrix $A$. If $p$ is a polynomial of degree $d$, then the sparsity pattern of $p(A)$ is contained in the union of the sparsity patterns of $I_n,A,\ldots,A^d$, so it corresponds to the vertex pairs in $G(A)$ connected by walks of length at most $d$. A slow growth of graph neighborhoods therefore leads to a slow fill-in for low-degree polynomials in $A$. Thus, the full-size global Lanczos approach is particularly suited to paths, low-dimensional grids, road networks, and sparse matrices from local PDE discretizations. In contrast, it is less suited to small-world graphs, such as social networks, where low powers of $A$ can already become nearly dense.

\begin{remark}
	\label{rem:full-size-global-lanczos-in-literature}
	The choice $B = I_n$ and its interpretation through orthogonal matrix polynomials is briefly discussed in \cite[Section~3.1]{FDFM23}. However, the computational benefits of exploiting sparsity in the basis matrices are not investigated there, and their numerical experiments only use smaller blocks, similar to \cite{ReichelGlobal15}. Here, we explicitly exploit the sparsity of the basis matrices $\vec V_j = p_{j-1}(A)$ to perform the full-size global Lanczos iterations efficiently via sparse matrix arithmetic.
\end{remark}

\subsection{Gauss--Radau error bounds}
\label{sec:GaussRadau}

We briefly recall how Gauss--Radau quadrature rules can be used to bound the error of Gauss quadrature, under suitable assumptions on the signs of the derivatives of $f$; see, e.g.,~\cite[Chapter~6]{GolubMeurantBook10}. This approach is applied to the global Lanczos approximation in~\cite{ReichelGlobal15}. A Gauss--Radau rule prescribes a fixed node $\zeta$, taken as either a lower or an upper bound on $\sigma(A)$ depending on the function $f$. 
More specifically, we consider an approximation
\begin{equation}
	\label{eqn:gauss-radau-approx}
    \mathcal{I}[f,\B] \approx \psi_{m+1}^{R}[f, \vec B] := \|\B\|_F^2 \e_1^T f(\widetilde{T}_{m+1})\e_1
\end{equation}
where $\widetilde{T}_{m+1}$ is an extension of the matrix $T_m$ generated by~\cref{alg:global-lanczos}
\begin{equation*}
    \widetilde{T}_{m+1}
    :=
    \begin{bmatrix}
        T_m & \beta_{m+1} \e_m\\
        \beta_{m+1} \e_m^T & \omega_{m+1}
    \end{bmatrix},
\end{equation*}
with eigenvalues $\widetilde{\theta}_1, \dots, \widetilde{\theta}_m, \zeta$. 
Note that $\widetilde{T}_{m+1}$ coincides with $T_{m+1}$ except for the entry at index $(m{+}1,m{+}1)$, which is replaced by $\omega_{m+1}$, chosen such that the prescribed value $\zeta$ is an eigenvalue of $\widetilde{T}_{m+1}$. Note that $\beta_{m+1}$ is already available after the $m$th step of~\cref{alg:global-lanczos}, so $\widetilde{T}_{m+1}$ can be formed without performing an additional Lanczos iteration.
As shown in \cite[eq.~(6.8)]{GolubMeurantBook10}, $\omega_{m+1}$ can be computed by solving an $m \times m$ tridiagonal linear system:
\begin{equation*}
    \omega_{m+1} = \zeta + \beta_{m+1}^2
       \e_m^T(T_m-\zeta I)^{-1}\e_m.
\end{equation*}
The errors of the approximations $\psi_m^{GL}[f, \vec B]$ and $\psi^R_{m+1}[f, \vec B]$ have the following representations.
\begin{theorem}{\cite[eq.~(1.4.14) and eq.~(1.4.23)]{GautschiBook04}}
    \label{thm:error-formulas-quadrature}
    Let $A\in\R^{n\times n}$ be symmetric with $\sigma(A)\subset[a,b]$, let $\B\in\R^{n\times \ell}$ and assume that $f \in C^{2m+1}([a, b])$. Assume that~\cref{alg:global-lanczos} does not break down before step $m$. Let $\theta_1, \dots, \theta_m$ be the eigenvalues of $T_m$, and $\widetilde{\theta}_1,\dots,\widetilde{\theta}_m, \zeta$ be the eigenvalues of $\widetilde{T}_{m+1}$, with $\zeta \in \{a, b\}$. Then there exist $\xi^G\in (a,b)$ and $\xi^R \in [a, b]$ such that
    \begin{align}
        \label{eq:error-formula-Gauss}
        \mathcal{I}[f,\B] - \psi^{GL}_m[f, \vec B] 
        &= 
        \frac{f^{(2m)}(\xi^G)}{(2m)!} 
        \int_a^b \prod_{k=1}^m(x-\theta_k)^2\,\de\mu(x) 
        \\
        \label{eq:error-formula-Gauss-Radau}
        \mathcal{I}[f,\B] - \psi^R_{m+1}[f, \vec B] 
        &=
        \frac{f^{(2m+1)}(\xi^R)}{(2m+1)!}
        \int_a^b (x- \zeta)\prod_{k=1}^m (x-\widetilde{\theta}_k)^2\,\de\mu(x). 
    \end{align} 
\end{theorem}

\cref{thm:error-formulas-quadrature} tells us whether $\psi^{GL}_m[f, \vec B]$ and $\psi^R_{m+1}[f, \vec B]$ bound $\mathcal{I}[f,\B]$ from above or below.
Assuming that
\begin{itemize}
	\item $f^{(2m)} \ge 0$ and $f^{(2m+1)} \ge 0$, with $\zeta = b$, or
	\item $f^{(2m)} \ge 0$ and $f^{(2m+1)} \le 0$, with $\zeta = a$,
\end{itemize}
we have 
\begin{equation*}
	\psi_m^{GL}[f, \vec B] \le \mathcal{I}[f, \vec B] \le \psi_m^R[f, \vec B].
\end{equation*}
Under these assumptions, we have the error bound 
\begin{equation}
	\label{eqn:gauss-gauss-radau-error-bound}
    \abs{\mathcal{I}[f,\B]-\psi_m^{GL}[f,\B]}
    \leq
    \abs{\psi_{m+1}^{R}[f,\B]-\psi_m^{GL}[f,\B]}.
\end{equation}
These conditions on the signs of derivatives are satisfied for several functions of interest, such as absolutely monotonic and completely monotonic functions; see, e.g., \cite[Chapter~IV]{WidderBook41}. 
We recall that a function $f$ is absolutely monotonic on the interval $[a, b]$ if
$f^{(k)}(x)\geq 0$ for all $x\in [a, b]$ and all $k\geq 0$, and completely
monotonic if $(-1)^k f^{(k)}(x)\geq 0$. For example, the exponential $f(x) = e^x$ is an absolutely monotonic function, and inverse fractional powers $f(x) = x^{-\alpha}$, $\alpha > 0$ are completely monotonic functions.

\subsection{Exact error formulation}
\label{subsec:exact-error-formulation}
 
The following elementary lemma shows that the error of the global Lanczos approximation \cref{eq:GL-approximation} to $\mathcal{I}[f, \vec B]$ can be written in terms of another trace.

\begin{lemma}
	\label{lemma:hermite-interpolant-quadrature}
	Let $q$ be a polynomial of degree at most $2m-1$ such that $q$ interpolates $f$ at all Ritz values $\theta_i$, i.e., $f(\theta_i)=q(\theta_i),\; 1 \leq i \leq m$. Then,
	\begin{align}
		\label{eq:classical-gauss-error}
		\mathcal{I}[f,\B] - \psi_m^{GL}[f, \vec B] &= \mathcal{I}[f-q,\B].
	\end{align}
\end{lemma}
\begin{proof}
	The fact that the global Lanczos approximation is a Gauss quadrature rule (Theorem~\ref{thm:gauss-approximation}) gives $\psi^{GL}_m[q, \vec B] = \mathcal{I}[q,\B]$. Hence,
	\begin{equation*}
		\mathcal{I}[f,\B]-\psi_m^{GL}[f, \vec B] = \mathcal{I}[f-q,\B] - \psi_m^{GL}[f-q, \vec B],
	\end{equation*}
	and the lemma is proved by noting that $q$ interpolates $f$ at the Ritz values.
\end{proof}
A natural candidate for this polynomial $q$ that uses all $2m$ degrees of freedom is the Hermite interpolant $h_m$ defined as
\begin{equation*}
	h_m(\theta_j) = f(\theta_j),
	\qquad
	h_m'(\theta_j) = f'(\theta_j),
	\qquad j = 1,\dots,m.
\end{equation*}

The next theorem gives a more explicit error formulation in terms of the spectral decomposition of the recurrence matrix $T_m$.
\begin{theorem}
\label{thm:errorformula}
Let $f \in C^{2}$ on an open interval containing $[\lambda_{\min}, \lambda_{\max}]$ and let $\eta_m(A) \coloneqq f(A) - h_m(A)$. Then,
\begin{equation}
	\label{eq:error}
\mathcal{I}[f,\boldsymbol{B}] - \psi_m^{GL}[f, \vec B]
  = \trace(\boldsymbol{B}^T \eta_m(A) \boldsymbol{B})
  = \|\boldsymbol{B}\|_F^2 \trace(\boldsymbol{V}_{m+1}^T g_m(A) \boldsymbol{V}_{m+1}),
\end{equation}
where $g_m$ is expressed in terms of divided differences \cite{deBoor05} of $f$,
\begin{equation*}
  g_m(z) := (\beta_2 \cdots \beta_{m+1})^2 f[\theta_1,\theta_1,\ldots,\theta_m,\theta_m,z].
\end{equation*}
\end{theorem}
\begin{proof}
The Newton form of the Hermite interpolation error is (see, e.g.,~\cite[Section~2.1.5]{StoerBulirsch02}) 
\begin{equation*}
	\eta_m(z)=f(z)-h_m(z)=\prod_{i=1}^m (z-\theta_i)^2 \, f[\theta_1,\theta_1,\ldots,\theta_m,\theta_m,z].
\end{equation*}
Lemma~\ref{lemma:hermite-interpolant-quadrature} with $q=h_m$ gives us 
\begin{equation*}
	\mathcal{I}[f,\boldsymbol{B}] - \psi_m^{GL}[f, \vec B] = \trace(\boldsymbol{B}^T\eta_m(A)\boldsymbol{B}) = \norm{\vec B}_F^2 \trace(\boldsymbol{V}_{m+1}^T g_m(A) \boldsymbol{V}_{m+1}),
\end{equation*}
where we used that $p_m(z) = (\beta_1\cdots\beta_{m+1})^{-1}\prod_{i=1}^m (z-\theta_i)$, $\boldsymbol{V}_{m+1} = p_m(A)\boldsymbol{B}$ and $\beta_1 = \|\boldsymbol{B}\|_F$.
\end{proof}

In many cases, $f$ can be written as a Cauchy-type integral or is closely related to one, see~\cref{subsubsec:evaluation-gm}, where we derive such representations for the functions used in our numerical experiments. By a Cauchy-type integral we mean that $f$ is analytic on a region $\Omega \subset \C$ and admits the integral representation
\begin{equation}\label{eq:CauchyStieltjes}
	f(z) = \int_{\Gamma}\frac{g(t)}{t-z}\,\mathrm{d}t, \qquad z\in \Omega,
\end{equation}
with a path $\Gamma \subset \C \setminus \Omega$ and function $g\colon \Gamma \to \C$. As in Frommer et al.~\cite{FrommerGuttelSchweitzer2014}, we mainly consider two cases of Cauchy-type representations:
\begin{enumerate}
	\item[(i)] $f$ is analytic on a region $\Omega' \supset \Omega$ and $\Gamma \subset \Omega'$ is a closed contour winding once around $\Omega$. By the Cauchy integral formula, \eqref{eq:CauchyStieltjes} then holds with $g(t) = \frac{1}{2\pi i} f(t)$.
	\item[(ii)] $f$ is a Stieltjes function, i.e.,
	\begin{equation}\label{eq:stieltjes}
		f(z) = \int_0^\infty \frac{\rho(t)}{z+t}\,\mathrm{d}t, \qquad z \in \Omega = \C \setminus (-\infty, 0],
	\end{equation}
	with $\rho$ of constant sign on $[0,\infty)$. After the substitution $t \mapsto -t$, this is~\eqref{eq:CauchyStieltjes} with $\Gamma = (-\infty, 0]$ and $g(t) = -\rho(-t)$.
\end{enumerate}
The next proposition shows that, for $f$ of the form~\eqref{eq:CauchyStieltjes}, the function $g_m$ also admits an integral representation.
\begin{proposition}\label{prop:errorformula-integral}
	Let $f$ be of the form~\eqref{eq:CauchyStieltjes} with $[\lambda_{\min}, \lambda_{\max}] \subset \Omega$. Then, for $z \in \Omega$,
\begin{equation}\label{eq:contourintegralerror}
	g_m(z) = (\beta_2 \cdots \beta_{m+1})^2 \int_\Gamma \frac{g(t)}{\det(tI - T_m)^2} \, \frac{\mathrm{d}t}{t - z}.
\end{equation}
\end{proposition}
\begin{proof}
	Since divided differences are linear, it suffices to compute the divided differences of $r_t(z) \coloneqq (t-z)^{-1}$ for fixed $t \in \Gamma$. A simple induction on the number of nodes shows that
\begin{equation*}
	r_t[x_0, \ldots, x_k] = \prod_{j=0}^k \frac{1}{t-x_j},
\end{equation*}
and by continuity this identity also holds for repeated nodes. Applying this to the definition of $g_m$, with nodes $\theta_1,\theta_1,\ldots,\theta_m,\theta_m \in [\lambda_{\min},\lambda_{\max}] \subset \Omega$, and using $\det(tI-T_m) = \prod_{j=1}^m(t-\theta_j)$ gives~\eqref{eq:contourintegralerror}.
\end{proof}

This error formulation also gives us a posteriori error bounds which can be computed by finding the extrema of $g_m$ over the spectral interval of $A$.

\begin{corollary}\label{cor:residuebounds}
	Under the same assumptions as Theorem~\ref{thm:errorformula}, we have
	\begin{equation}
        \label{eq:aposteriori-bound}
        \norm{\B}_F^2 \min_{z\in[\lambda_{\min},\lambda_{\max}]} g_m(z) \;\leq\; \mathcal{I}[f,\B]- \psi^{GL}_m[f, \vec B] \;\leq\; \norm{\B}_F^2 \max_{z\in[\lambda_{\min},\lambda_{\max}]} g_m(z).
    \end{equation}
\end{corollary}
\begin{proof}
	Let $A = \sum_{i=1}^n \lambda_i \vec z_i \vec z_i^T$ be a spectral decomposition of $A$ with orthonormal eigenvectors $\vec z_i$, so that $g_m(A) = \sum_{i=1}^n g_m(\lambda_i)\vec z_i \vec z_i^T$. Substituting into~\cref{eq:error} yields
\begin{equation*}
	\mathcal{I}[f,\B] - \psi^{GL}_m[f, \vec B] = \norm{\B}_F^2 \sum_{i=1}^n g_m(\lambda_i)\, w_i,
\end{equation*}
where $w_i := \trace(\vec V_{m+1}^T \vec z_i \vec z_i^T \vec V_{m+1}) = \norm{\vec V_{m+1}^T \vec z_i}_2^2 \ge 0.$
The weights are nonnegative and sum to one, since $\sum_{i=1}^n \vec z_i \vec z_i^T = I_n$ gives $\sum_{i=1}^n w_i = \trace(\vec V_{m+1}^T \vec V_{m+1}) = \norm{\vec V_{m+1}}_F^2 = 1.$
Hence, $\big(\mathcal{I}[f,\B] - \psi^{GL}_m\big)/\norm{\B}_F^2$ is a convex combination of the values $g_m(\lambda_i)$, with all nodes $\lambda_i \in \sigma(A)$, and thus lies between the minimum and maximum of $g_m$.
\end{proof}
\begin{remark}\label{rem:gauss-radau-vs-residue-bounds}
    Under the derivative sign conditions stated after \cref{thm:error-formulas-quadrature}, the Gauss--Radau error bound~\cref{eqn:gauss-gauss-radau-error-bound} with $[a, b] = [\lambda_{\min}, \lambda_{\max}]$  is never worse than the upper bound in~\cref{cor:residuebounds}.
    Indeed, let $\xi_j\in[a,b]$ and $\widetilde w_j>0$, $j=1,\dots,m+1$, be the nodes and weights of the Gauss--Radau rule, so that $\psi_{m+1}^R[u,\B]=\sum_{j=1}^{m+1}\widetilde w_j u(\xi_j)$. Since this rule is exact for polynomials of degree at most $2m$, we have $\psi_{m+1}^R[h_m,\B]=\mathcal I[h_m,\B]=\psi_m^{GL}[f,\B]$ and hence $\psi_{m+1}^R[f, \vec B]  - \psi_{m}^{GL}[f, \vec B] = \psi_{m+1}^R[f - h_m, \vec B]$. Using $f-h_m=\norm{\B}_F^2p_m^2g_m$ from the proof of~\cref{thm:errorformula}, we obtain
    \begin{equation}
		\label{eqn:rem-gauss-radau-vs-residue-bounds}
        \psi_{m+1}^R[f,\B]-\psi_m^{GL}[f,\B]
        =\norm{\B}_F^2\sum_{j=1}^{m+1}\widetilde w_j p_m(\xi_j)^2g_m(\xi_j),
    \end{equation}
    and because of polynomial exactness we have
    \begin{equation*}
        \sum_{j=1}^{m+1}\widetilde w_j p_m(\xi_j)^2
        =\int_a^b p_m(x)^2\,\de\mu(x)=1.
    \end{equation*}
    So the sum in the right-hand side of \cref{eqn:rem-gauss-radau-vs-residue-bounds} is a convex combination of values of $g_m$ on $[a, b]$, and we have
    \begin{align*}
        \mathcal I[f,\B]-\psi_m^{GL}[f,\B]
        \leq\psi_{m+1}^R[f,\B]-\psi_m^{GL}[f,\B]
        \leq\norm{\B}_F^2\max_{x\in[a,b]}g_m(x).
    \end{align*}
\end{remark}

\section{Other trace estimators}
\label{sec:other-trace-estimators}

In this section we provide a brief overview of methods for estimating the trace of $f(A)$, focusing on stochastic trace estimators \cite{Hutchinson89,EpperlyTroppWebber24,MeyerMuscoMuscoWoodruff2021} and probing methods \cite{FrommerSchimmelSchweitzer21,FrommerRinelliSchweitzer25}. 
Several trace estimators have the form 
\begin{equation}
    \label{eq:trace-estimation-quadforms}
    \trace(f(A)) \approx \sum_{i=1}^s \vec x_i^T f(A) \vec x_i = \trace(\vec X^T f(A) \vec X) = \mathcal{I}[f, \vec X],
\end{equation}
for suitable vectors $\vec x_1, \dots, \vec x_s$, with $\vec X = [\vec x_1, \dots, \vec x_s]$. The right-hand side of \cref{eq:trace-estimation-quadforms} can then be approximated with \cref{eq:GL-approximation}. Hence, approximations of $\trace(f(A))$ via~\cref{eq:trace-estimation-quadforms} necessarily go on two levels, which provide two different sources for the error: the approximation with the exact trace estimator, and the approximation of the quadratic forms with the Lanczos method. Some estimators~\cite{EpperlyTroppWebber24,MeyerMuscoMuscoWoodruff2021} also require the approximation of matrix-vector multiplications with $f(A)$ inside the process.

\subsection{Hutchinson's estimator}
\label{subsec:hutchinson}

The simplest stochastic estimator to estimate the trace of a matrix $M$ is Hutchinson's estimator~\cite{Hutchinson89,Girard1989}. The estimator is the following:
\begin{equation}
	\label{eq:Hutchinson-estimator}
    \trace(M) \approx \frac{1}{s} \sum_{i=1}^s \vec x_i^T M \vec x_i=:\trace^H_s(M),
\end{equation}
where $\vec x_i$ are i.i.d. random vectors that satisfy $\mathbb{E} [\vec x_i] = \vec 0$ and $\mathbb{E}[\vec x_i\vec x_i^T] = I_n$.
Popular choices are Gaussian vectors (i.i.d. $\mathcal{N}(0,1)$ entries) and Rademacher vectors (i.i.d. $\pm1$ entries). The estimator is unbiased since, under the above assumptions, $\mathbb{E}[\vec x_i^T M \vec x_i] = \trace(M)$.
The quality of unbiased stochastic estimators can be analyzed via their variance, which for Hutchinson's estimator is given by
\begin{equation} \label{eq:variance-hutchinson}
    \mathbb{V}[\trace^H_s(M)] =
    \begin{cases}
        \frac{2}{s}\|M\|_F^2 & \text{for Gaussian vectors},\\
        \frac{2}{s}\|M-\operatorname{diag}(M)\|_F^2 & \text{for Rademacher vectors}.
    \end{cases}
\end{equation}
In order to have $\mathbb P(|\trace(M)-\trace^H_s(M)|\geq \epsilon)\leq \delta$ for a target accuracy $\epsilon$ and prescribed failure probability $\delta$, we need to choose
\begin{equation} \label{eq:tailbound-hutchinson}
    s = O\left( \frac{\log(1/\delta)}{\epsilon^2}\|M\|_F^2 \right),
\end{equation}
see, e.g.,~\cite{CortinovisKressner22}.

\subsection{XTrace and XNysTrace}
\label{subsec:xtrace}
Epperly, Tropp and Webber~\cite{EpperlyTroppWebber24} designed \texttt{XTrace} and \texttt{XNysTrace}, two methods which are based on the same principle as the \texttt{Hutch++} estimator~\cite{MeyerMuscoMuscoWoodruff2021}. The idea is to compute the trace of a randomized low-rank approximation $\vec Q \vec Q^T M$ of $M$ exactly, and use Hutchinson's estimator to approximate the remainder term $(I - \vec Q\vec Q^T)M(I - \vec Q\vec Q^T)$. The variance, and thus error, decreases considerably whenever $M$ has fast eigenvalue decay. In contrast to \texttt{Hutch++}, \texttt{XTrace} and \texttt{XNysTrace} satisfy the \emph{exchangeability principle}, meaning that the estimator of $\trace(M)$ is a symmetric function of the i.i.d test vectors $\vec x_1, \dots,\vec x_s$.

Let $\vec X = [\vec x_1,\dots,\vec x_s]$, and let $\vec X_{-i}$ denote $\vec X$ with its $i$th column removed. \texttt{XTrace} sets $\vec Q_{(i)} = \orth(M\vec X_{-i})$ and defines the basic estimators
\begin{equation}
	\label{eq:xtrace-basic}
    t_i := \trace(\vec Q_{(i)}^T M \vec Q_{(i)}) + \vec x_i^T (I - \vec Q_{(i)}\vec Q_{(i)}^T) M (I - \vec Q_{(i)}\vec Q_{(i)}^T) \vec x_i, \qquad i = 1,\dots,s.
\end{equation}
Since $\vec x_i$ is independent of $\vec Q_{(i)}$, each $t_i$ is unbiased, and so is their average
\begin{equation}
	\label{eq:xtrace-estimator}
    \trace^{X}_s(M) := \frac{1}{s}\sum_{i=1}^s t_i,
\end{equation}
which is a symmetric function of $\vec x_1,\dots,\vec x_s$. The authors of~\cite{EpperlyTroppWebber24} provide an efficient implementation that only needs two products of $M$ with $n\times s$ block vectors. It first forms $\vec Y = M\vec X$ with thin QR factorization $\vec Y = \vec Q\vec R$, and then $\vec Z = M\vec Q$. Every $\vec Q_{(i)}$ is obtained from $\vec Q$ and $\vec R$ by a rank-one downdate, so all basic estimators $t_1,\dots,t_s$ can be read off from $\vec X$, $\vec Q$, $\vec R$ and $\vec Z$ with $O(s^2 n)$ additional operations. This amounts to $2s$ matrix-vector products with $M$, the same asymptotic cost as \texttt{Hutch++}.

When $M$ is symmetric positive semidefinite (e.g., $M = f(A)$ with $f \geq 0$ on $\sigma(A)$), \texttt{XNysTrace} replaces the randomized SVD by the Nystr\"om approximation $M\langle \vec X\rangle := (M\vec X)(\vec X^T M\vec X)^{\dagger}(M\vec X)^T$, which needs only the single block product $M\vec X$. The basic estimators and the \texttt{XNysTrace} estimator are
\begin{equation}
	\label{eq:xnystrace-estimator}
    t_i := \trace(M\langle \vec X_{-i}\rangle) + \vec x_i^T\bigl(M - M\langle \vec X_{-i}\rangle\bigr)\vec x_i, \qquad \trace^{XN}_s(M) := \frac{1}{s}\sum_{i=1}^s t_i.
\end{equation}
For a budget of $m$ matrix-vector products, \texttt{Hutch++}, \texttt{XTrace} and \texttt{XNysTrace} build low-rank approximations of rank about $m/3$, $m/2$ and $m$, respectively. Suppose the eigenvalues of $M$ decay like $\alpha^i$, their errors decay (up to polynomial factors in $m$) like $\alpha^{m/3}$, $\alpha^{m/2}$ and $\alpha^{m}$~\cite{EpperlyTroppWebber24}. Moreover, the spread of the basic estimators $t_i$ provides the a posteriori error estimate
\begin{equation}
	\label{eq:xtrace-error-estimate}
    \widehat{\epsilon}^{\,2} := \frac{1}{s(s-1)}\sum_{i=1}^s \bigl(t_i - \trace^{X}_s(M)\bigr)^2 \approx \bigl|\trace(M) - \trace^{X}_s(M)\bigr|^2,
\end{equation}
and analogously for \texttt{XNysTrace}. The estimate is typically a modest underestimation of the error, because the $t_i$ are often weakly positively correlated, as observed in \cite{EpperlyTroppWebber24}. In our setting $M = f(A)$, so each block product with $M$ must itself be approximated, which we do with the Lanczos method applied to all columns of the block simultaneously. Hence \texttt{XTrace} requires the two approximate block products $f(A)\vec X$ and $f(A)\vec Q$, whereas \texttt{XNysTrace} requires only $f(A)\vec X$. Approximating these products much more accurately than the inherent sampling error of the stochastic estimator wastes computational effort, while approximating them much less accurately makes the Lanczos error dominate the overall error. We therefore discuss practical stopping criteria for these inner Lanczos iterations in~\cref{subsubsec:inner-stopping}. 

\subsection{Probing methods}
\label{subsec:probing-methods}

When $A$ is a sparse matrix, the trace of $f(A)$ can be approximated via probing methods \cite{FrommerSchimmelSchweitzer21}, which construct probing vectors $\x_i$ based on a distance-$d$ coloring associated with $A$. The estimator has the form
\begin{equation}
    \label{eq:deterministic-probing}
    \mathcal{T}_d(f(A))
    :=\sum_{i=1}^{s} \x_s^T f(A)\x_s
    =\mathcal{I}[f,\X],
\end{equation}
where $\X:=[\x_1,\ldots,\x_{s}]$. One can then use global Lanczos or several single-vector Lanczos runs to get the approximation of the trace.

Probing is effective when the entries of $f(A)$ decay with graph distance. We refer to
\cite{FrommerSchimmelSchweitzer21} for an analysis of deterministic probing for
matrix functions.

In deterministic probing, the entries of the vectors are only $0$ and $1$.
Stochastic probing~\cite{FrommerRinelliSchweitzer25} replaces the nonzero entries by independent Rademacher variables which take uniformly the values $\pm 1$. This typically leads to a smaller error than deterministic probing; see
\cite{FrommerRinelliSchweitzer25} for a detailed analysis.

For a general graph, computing a coloring with a minimal number of colors is already NP-hard for $d=1$. In practice, distance-$d$ colorings can be constructed using greedy procedures for general graphs, while explicit constructions are available in special cases, such as paths and regular grids; see~\cite[Section~2]{FrommerSchimmelSchweitzer21}.

\section{Stochastic error correction}
\label{sec:stochastic-error-correction}

In this section we propose new hybrid trace estimators that combine the deterministic full-size ($\B=I_n$) global Lanczos approximation $\psi_m^{GL}$~\cref{eq:full-size-GL} with a stochastic correction of its error.
The key observation is that, by~\cref{eq:error}, the error of this approximation is itself the trace of a matrix function, and can therefore be estimated by any stochastic trace estimator. Since this matrix is typically much smaller in Frobenius norm than $f(A)$, the correction has a considerably lower variance than a stochastic estimator applied directly to $\trace(f(A))$.
Each choice of stochastic estimator gives rise to a different hybrid method; we focus on \texttt{GLHutch}, \texttt{GLXTrace}, and \texttt{GLXNysTrace} which are based on Hutchinson's estimator~\cite{Hutchinson89} and on \texttt{XTrace}/\texttt{XNysTrace}~\cite{EpperlyTroppWebber24}, respectively.
We first introduce these methods, and then discuss their theoretical properties and implementation aspects.

\subsection{Methods}
\label{subsec:methods}

Let $M \coloneqq \eta_m(A) = f(A) - h_m(A)$ denote the \emph{error matrix} of the $m$-th full-size global Lanczos approximation. Since $\B = I_n$, we have $\norm{\B}_F^2 = n$ and $\V_{m+1} = p_m(A)$, so that by~\cref{eq:error}, we have $M = n\, p_m(A)^2 g_m(A)$ and
\begin{equation}
\label{eqn:error-of-global-lanczos-as-trace}
\trace(f(A)) - \psi_m^{GL} = \trace(M).
\end{equation}
\subsubsection*{GLHutch}
Applying Hutchinson's estimator~\cref{eq:Hutchinson-estimator} to $M$ gives the \emph{corrected estimate}
\begin{equation}
	\label{eq:corrected-estimate}
	\trace(f(A)) \approx \psi_m^{GL} + \trace^H_s(M)
	= \psi_m^{GL} + \frac{1}{s}\sum_{i=1}^s \x_i^T M \x_i,
\end{equation}
with i.i.d.~random vectors $\x_1,\dots,\x_s \in \R^n$ as in~\cref{subsec:hutchinson}. Letting $\X = [\x_1, \ldots, \x_s] \in \R^{n \times s}$, the correction term in~\cref{eq:corrected-estimate} is $\trace^H_s(M) = \frac{1}{s}\trace(\X^T M \X)$. Using $M = n\,\V_{m+1}^T\, g_m(A)\, \V_{m+1}$, we approximate this term using a second run of global Lanczos~\cref{eq:GL-approximation}
\begin{equation}
	\label{eq:correction-as-quadratic-form}
	\trace_s^H(M)
	= \frac{n}{s} \trace\big((\V_{m+1}\X)^T g_m(A)\, (\V_{m+1}\X)\big)
	\approx \frac{n}{s}\, \psi^{GL}_r[g_m, \V_{m+1}\X],
\end{equation}
where the second run is started from the block $\V_{m+1}\X$ and the number of iterations $r$ is determined by an inner stopping criterion, discussed in~\cref{subsubsec:inner-stopping}. The starting block is cheap to form: since $\V_{m+1} =p_m(A)$, it can be computed with that same recurrence and the known coefficients $\alpha_k, \beta_{k+1}$. This requires only $m$ sparse products with $A$ and avoids directly computing the product with $\V_{m+1}$.

\subsubsection*{GLXTrace and GLXNysTrace}
If the error matrix $M$ is numerically low-rank, then we expect a considerable variance reduction by applying \texttt{XTrace} or \texttt{XNysTrace} to $M$ compared to using Hutchinson.
Applying \texttt{XTrace}~\cref{eq:xtrace-estimator} to $M$ gives the corrected estimate
\begin{equation}
	\label{eq:corrected-estimate-xtrace}
	\trace(f(A)) \approx \psi_m^{GL} + \trace^{X}_s(M).
\end{equation}
Recall from~\cref{subsec:xtrace} that \texttt{XTrace} accesses $M$ only through the two block products $M\X$ and $M\vec Q$, with $\vec Q = \orth(M\X) \in \R^{n\times s}$. For any block vector $\W \in \R^{n\times s}$ we compute
\begin{equation*}
	M\W = n\, p_m(A)\, g_m(A)\, p_m(A)\W,
\end{equation*}
performing operations from right to left: we form $p_m(A)\W$ with the three-term recurrence as above, approximate the product with $g_m(A)$ with the Lanczos method, and apply $p_m(A)$ once more. The Lanczos method is run on the $s$ columns simultaneously, so that each iteration needs a single product of $A$ with an $n\times s$ block. Each block product with $M$ therefore costs $2m$ sparse products of $A$ with an $n\times s$ block, plus the $r$ such products of the Lanczos run for $g_m(A)$, where $r$ is set by the inner stopping criterion of~\cref{subsubsec:inner-stopping}.

\texttt{XNysTrace} additionally requires $M$ to be semidefinite. The next proposition shows that this holds under the same sign condition on $f^{(2m)}$ as the Gauss--Radau bounds in~\cref{sec:GaussRadau}.

\begin{proposition}
	\label{prop:M-semidefinite}
	Let $f \in C^{2m}$ on an open interval containing $[\lambda_{\min},\lambda_{\max}]$, with $f^{(2m)} \ge 0$ (resp.\ $f^{(2m)} \le 0$) on $[\lambda_{\min},\lambda_{\max}]$. Then the error matrix $M$ is symmetric positive (resp.\ negative) semidefinite.
\end{proposition}
\begin{proof}
	The Ritz values $\theta_1,\dots,\theta_m$ lie in $[\lambda_{\min},\lambda_{\max}]$. By the mean value theorem for divided differences~\cite{deBoor05}, for every $z \in [\lambda_{\min},\lambda_{\max}]$ there exists $\xi$ in this interval such that $f[\theta_1,\theta_1,\ldots,\theta_m,\theta_m,z] = f^{(2m)}(\xi)/(2m)!$. Hence $g_m(z)$ has the sign of $f^{(2m)}$, and so do the eigenvalues $n\,p_m(\lambda_i)^2 g_m(\lambda_i)$, $i = 1,\dots,n$, of $M = n\,p_m(A)^2 g_m(A)$.
\end{proof}

For example, $f(x) = e^x$ and $f(x) = x^{-\alpha}$ give a positive semidefinite $M$, while $f(x) = \log(x)$ gives a negative semidefinite $M$. Applying \texttt{XNysTrace}~\cref{eq:xnystrace-estimator} to $M$, or to $-M$ in the negative semidefinite case, gives the corrected estimate
\begin{equation}
	\label{eq:corrected-estimate-xnystrace}
	\trace(f(A)) \approx \psi_m^{GL} + \trace^{XN}_s(M)
	\qquad\text{or}\qquad
	\trace(f(A)) \approx \psi_m^{GL} - \trace^{XN}_s(-M),
\end{equation}
respectively. \texttt{XNysTrace} only requires one block product $M\X$ instead of two for \texttt{XTrace}. \Cref{alg:gl-correction} summarizes the three hybrid methods.

\begin{algorithm}[t]
\caption{Global Lanczos with stochastic error correction}
\label{alg:gl-correction}
\begin{algorithmic}[1]
\REQUIRE Symmetric $A\in\R^{n\times n}$, function $f$, GL steps $m$, number of test vectors $s$, estimator $\mathrm{E}\in\{\trace^H_s, \trace^X_s, \trace^{XN}_s\}$, initial inner tolerance $\tau_0$, factor $c \in (0,1)$
\ENSURE Approximation of $\trace(f(A))$
\STATE $(\psi_m^{GL}, T_m) \leftarrow \texttt{GlobalLanczos}(A, I_n, f, m)$ \COMMENT{\cref{alg:global-lanczos}}
\STATE Draw i.i.d.~Gaussian or Rademacher $\X = [\x_1,\dots,\x_s]\in\R^{n\times s}$
\IF{$\mathrm{E} = \trace_s^H$}
	\STATE $\vec P = p_m(A)\X$
	\STATE $\psi_r^{GL}[g_m,\vec P] \leftarrow \texttt{GlobalLanczos}(A, \vec P, g_m, r)$ \COMMENT{$r$ set by~\cref{eq:inner-stopping-hutch}} \label{line:hutch-gl}
	\RETURN $\psi_m^{GL} + \frac{n}{s}\psi_r^{GL}[g_m,\vec P]$
\ELSIF{$\mathrm{E} = \trace_s^X$}
	\STATE $\tau = \tau_0$
	\STATE $\Y \approx M\X$, thin QR $\Y = \vec Q\vec R$, $\vec Z \approx M\vec Q$, both with tolerance $\tau$ \label{line:xtrace-products}
	\STATE Compute $t_1,\dots,t_s$~\cref{eq:xtrace-basic} from $\X,\vec Q,\vec R,\vec Z$ and $\widehat\epsilon$~\cref{eq:xtrace-error-estimate} \label{line:xtrace-estimate}
	\STATE Compute the correction estimate $\widehat{\mathrm{tr}} = s^{-1} \sum_{i=1}^s t_i$ 
	\IF{$c\,\widehat\epsilon/\abs{\widehat{\mathrm{tr}}} < \tau_0$} \label{line:xtrace-check}
		\STATE $\tau = c\,\widehat\epsilon/\abs{\widehat{\mathrm{tr}}}$ and repeat lines~\ref{line:xtrace-products}--\ref{line:xtrace-estimate} \label{line:xtrace-refine}
	\ENDIF
	\RETURN $\psi_m^{GL} + \frac{1}{s}\sum_{i=1}^s t_i$, error estimate $\widehat\epsilon$
\ELSIF{$\mathrm{E} = \trace_s^{XN}$}
	\STATE $\sigma = \pm 1$, the sign of $f^{(2m)}$ on $[\lambda_{\min},\lambda_{\max}]$ \COMMENT{$\sigma M$ is PSD by \cref{prop:M-semidefinite}}
	\STATE $\Y \approx \sigma M\X$ with tolerance $\tau_0$ \label{line:xnys-product}
	\STATE Compute $t_1,\dots,t_s$~\cref{eq:xnystrace-estimator} for $\sigma M$ from $\X,\Y$ and $\widehat\epsilon$ \label{line:xnys-estimate}
	\IF{$c\,\widehat\epsilon/\abs{\widehat{\mathrm{tr}}} < \tau_0$} \label{line:xnys-check}
		\STATE Continue line~\ref{line:xnys-product} to tolerance $c\,\widehat\epsilon/\abs{\widehat{\mathrm{tr}}}$ and repeat line~\ref{line:xnys-estimate} \label{line:xnys-refine}
	\ENDIF
	\RETURN $\psi_m^{GL} + \frac{\sigma}{s}\sum_{i=1}^s t_i$, error estimate $\widehat\epsilon$
\ENDIF
\end{algorithmic}
\end{algorithm}

\subsection{Theoretical aspects}

By~\cref{lemma:hermite-interpolant-quadrature}, the error~\cref{eqn:error-of-global-lanczos-as-trace} equals $\trace(f(A)-q(A))$ for any $q\in\poly_{2m-1}$ that interpolates $f$ at the Ritz values $\theta_1,\dots,\theta_m$, not only for $q = h_m$. In particular, the Lagrange interpolant $\ell_m\in\poly_{m-1}$ of $f$ at $\theta_1,\dots,\theta_m$ gives an error matrix $M_\ell \coloneqq f(A)-\ell_m(A)$ with $\trace(M_\ell) = \trace(M)$, so the stochastic correction could equally be applied to $M_\ell$. We now argue why $M$ is the better choice.

Since the variance of Hutchinson's estimator with Gaussian vectors is proportional to the squared Frobenius norm, see~\cref{eq:variance-hutchinson}, we compare $\norm{M}_F$ and $\norm{M_\ell}_F$. For the Hermite choice, the next proposition provides an upper bound.

\begin{proposition}\label{prop:variance-error-operator}
Let $a < b$ be such that $\sigma(A) \subset [a,b]$, and let $f \in
C^{2m}$ on an open interval containing $[a,b]$ with $f^{(2m)}$ of
constant sign there. Then
\begin{equation*}
\norm{M}_F^2 \;\le\; 2 n  \, \max_{x \in [a,b]} \abs{\eta_m(x)}
\cdot  E_{2m-1}(f,[a,b]).
\end{equation*}
\end{proposition}

\begin{proof}
By~\cref{prop:M-semidefinite}, $M = \eta_m(A)$ is semidefinite, so its eigenvalues $\eta_m(\lambda_1), \dots, \eta_m(\lambda_n)$ share the same sign and $\abs{\trace(M)} = \sum_{i} \abs{\eta_m(\lambda_i)}$. Hence
\begin{align*}
	\norm{M}_F^2 &\le \norm{M}_2 \, \abs{\trace (M)} \\
	&\le \, \max_{x \in [a, b]} \abs{\eta_m(x)} \cdot \abs*{\mathcal{I}[f, I_n] - \psi_m^{GL}[f, I_n]}.
\end{align*}
The statement then follows by using \cref{cor:convergence-Gauss-quadrature} to bound the right-hand side. 
\end{proof}

Both factors in the bound of~\cref{prop:variance-error-operator} are errors of polynomial approximations of degree $2m-1$, since $h_m\in\poly_{2m-1}$. For the Lagrange interpolant, $\ell_m\in\poly_{m-1}$ gives
\begin{equation}
	\label{eq:frobenius-lagrange}
	\norm{M_\ell}_F^2 = \sum_{i=1}^n (f-\ell_m)(\lambda_i)^2 \ge n\,\delta_{m-1}^2,
	\qquad
	\delta_{m-1} \coloneqq \min_{q\in\poly_{m-1}} \Big(\frac1n\sum_{i=1}^n (f-q)(\lambda_i)^2\Big)^{1/2},
\end{equation}
which only involves an approximation of degree $m-1$. This difference can be understood from the Newton forms of the two interpolation errors. With $\omega(z) = \prod_{j=1}^m (z-\theta_j)$, the Hermite error $\eta_m(z) = \omega(z)^2 f[\theta_1,\theta_1,\dots,\theta_m,\theta_m,z]$ does not change sign on $[a,b]$, whereas the Lagrange error $(f-\ell_m)(z) = \omega(z)\, f[\theta_1,\dots,\theta_m,z]$ may change sign at every Ritz value. Hence $M_\ell$ can have eigenvalues of both signs, which cancel in the trace but not in the Frobenius norm.

To see how the bounds of~\cref{prop:variance-error-operator} and~\cref{eq:frobenius-lagrange} behave as $m$ grows, let $f$ be analytic. Write $x = \frac{a+b}{2} + \frac{b-a}{2}t$ and expand $f(x) = \sum_{k\ge 0}\gamma_k T_k(t)$ in Chebyshev polynomials $T_k$. If $f$ is analytic and bounded by $C$ inside the Bernstein ellipse $E_\rho$, $\rho>1$, i.e., the ellipse with foci $a,b$ and semi-axes summing to $\rho\,\frac{b-a}{2}$, then $\abs{\gamma_k}\le 2C\rho^{-k}$~\cite[Theorem~8.1]{Trefethen19}, and therefore
\begin{equation*}
	E_{2m-1}(f,[a,b]) \le \sum_{k\ge 2m}\abs{\gamma_k} \le \frac{2C}{\rho-1}\,\rho^{-2m}.
\end{equation*}
By the Hermite integral formula for the interpolation error~\cite[Chapter~11]{Trefethen19}, we expect $\max_{[a,b]}\abs{\eta_m}$ to decay at the same rate $\rho^{-2m}$ when the Ritz values are well spread over $[a,b]$, so that~\cref{prop:variance-error-operator} gives $\norm{M}_F \lesssim \sqrt{n}\,\rho^{-2m}$. For~\cref{eq:frobenius-lagrange}, suppose that the eigenvalues fill $[a,b]$, in the sense that they are asymptotically distributed according to the Chebyshev density $1/(\pi\sqrt{1-t^2})$, as is the case, e.g., for the 1D discrete Laplacian. Then the $T_k$ are approximately orthogonal with respect to the eigenvalue distribution, and for $m\ll n$
\begin{equation*}
	\delta_{m-1}^2 \approx \frac12\sum_{k\ge m}\gamma_k^2 \ge \frac{\gamma_m^2}{2}.
\end{equation*}
Combining both bounds, we obtain, up to factors independent of $m$ and $n$,
\begin{equation*}
	\frac{\norm{M_\ell}_F}{\norm{M}_F} \gtrsim \abs{\gamma_m}\,\rho^{2m}.
\end{equation*}
If $\rho$ is taken as large as possible, the Chebyshev coefficients typically decay no faster than $\abs{\gamma_m}\approx\rho^{-m}$, so this ratio grows like $\rho^m$ and the number of samples needed for a given accuracy is larger by a factor of order $\rho^{2m}$ for $M_\ell$ than for $M$. For example, for $f(x) = \log x$ with $0<a<b$, the singularity at $0$ gives $\rho = (\sqrt\kappa+1)/(\sqrt\kappa-1)$ with $\kappa = b/a$, while for the entire function $f(x) = e^x$ the ratio grows faster than any geometric rate. In other words, the Lagrange error matrix sees the Chebyshev coefficient $\gamma_m$, whereas its trace only depends on $\gamma_{2m}, \gamma_{2m+1}, \ldots$, and the gap widens as $m$ increases.

Besides the smaller variance, the Hermite choice gives $M$ a structure that the methods of~\cref{subsec:methods} rely on. First, the factorization $M = n\,p_m(A)^2 g_m(A)$ writes the correction as a quadratic form in the block $\V_{m+1}\X = p_m(A)\X$, which allows us to compute it with a second run of global Lanczos~\cref{eq:correction-as-quadratic-form}. Second, by~\cref{prop:M-semidefinite}, $M$ is semidefinite whenever $f^{(2m)}$ has constant sign, which is required for \texttt{GLXNysTrace} and in~\cref{prop:variance-error-operator}. Neither of these properties holds for $M_\ell$: since $f-\ell_m$ contains the factor $\omega$ only once, $M_\ell$ admits no such symmetric factorization and is in general indefinite, so that \texttt{XNysTrace} cannot be applied.

\subsection{Implementation aspects}
\label{subsec:implementation-aspects}

The methods in~\cref{alg:gl-correction} require inner Lanczos runs for the function $g_m$: a single global Lanczos run for \texttt{GLHutch}, and a Lanczos run for each block product with $M$ for \texttt{GLXTrace} and \texttt{GLXNysTrace}. We first discuss when to stop these inner runs, and then how to evaluate $g_m$ accurately.

\subsubsection{Inner stopping criteria}
\label{subsubsec:inner-stopping}
The inner Lanczos runs are stopped differently for the Hutchinson correction and for the \texttt{XTrace}/\texttt{XNysTrace} corrections.

\emph{GLHutch.} The correction~\cref{eq:correction-as-quadratic-form} is computed with a single global Lanczos run on $(A,\vec P)$, $\vec P = p_m(A)\X$, see line~\ref{line:hutch-gl} of~\cref{alg:gl-correction}. Its error can therefore be estimated with the Gauss--Radau bound~\cref{eqn:gauss-gauss-radau-error-bound} applied to $g_m$. By~\cref{eq:tailbound-hutchinson}, the sampling error of Hutchinson's estimator decays like $1/\sqrt{s}$, so it suffices that the inner error is a fraction of this. We stop at the first $r$ for which
\begin{equation}
	\label{eq:inner-stopping-hutch}
	\frac{\abs{\psi^R_{r+1}[g_m,\vec P] - \psi_r^{GL}[g_m,\vec P]}}{\abs{\psi_r^{GL}[g_m,\vec P]}} \le \frac{C}{\sqrt{s}},
\end{equation}
with $C = 10^{-2}$ in our experiments.

\emph{GLXTrace/GLXNysTrace.} Each product with $M$ in lines~\ref{line:xtrace-products} and~\ref{line:xnys-product} of~\cref{alg:gl-correction} requires an approximation of $g_m(A)\vec P$, with $\vec P = [\vec p_1,\dots,\vec p_s] = p_m(A)\W$. For each column we run the Lanczos method on $(A,\vec p_i)$, which after $j$ iterations gives $g_m(A)\vec p_i \approx \vec y_j^{(i)} \coloneqq \norm{\vec p_i}_2\, U_j^{(i)} g_m(T_j^{(i)})\,\e_1$. The iteration with \emph{tolerance} $\tau$ stops at the first $j$ for which
\begin{equation}
	\label{eq:inner-stopping}
	\max_{i=1,\dots,s} \frac{\norm{\vec y_j^{(i)} - \vec y_{j-1}^{(i)}}_2}{\norm{\vec y_j^{(i)}}_2} \le \tau,
\end{equation}
i.e., when the relative change of every column between two consecutive iterations is at most $\tau$. Since the Krylov bases are nested and orthonormal, these norms can be computed from the coefficient vectors $g_m(T_j^{(i)})\e_1$ of length $j$, so checking~\cref{eq:inner-stopping} requires no operations with vectors of length $n$.

The tolerance $\tau$ should be chosen such that the error of the inner Lanczos approximations is small compared to the sampling error of the stochastic estimator: a much smaller $\tau$ wastes iterations, while a larger $\tau$ makes the inner error dominate the total error. Since the sampling error is not known in advance, we determine $\tau$ adaptively using the error estimate $\widehat\epsilon$~\cref{eq:xtrace-error-estimate}. First, an estimate $\widehat{\mathrm{tr}}$ is computed with a moderate tolerance $\tau_0$ (lines~\ref{line:xtrace-products}--\ref{line:xtrace-estimate} and~\ref{line:xnys-product}--\ref{line:xnys-estimate}). Heuristically, relative errors of size $\tau$ in the products with $M$ perturb $\widehat{\mathrm{tr}}$ by about $\tau\abs{\widehat{\mathrm{tr}}}$, so we require $\tau\abs{\widehat{\mathrm{tr}}} \le c\,\widehat\epsilon$ for a factor $c<1$. If $\tau_0$ already satisfies this, the first estimate is accepted. Otherwise, the products are refined with the tolerance $\tau = c\,\widehat\epsilon/\abs{\widehat{\mathrm{tr}}}$ (lines~\ref{line:xtrace-check}--\ref{line:xtrace-refine} and~\ref{line:xnys-check}--\ref{line:xnys-refine}). For \texttt{XTrace}, the second product $M\vec Q$ depends on the first through $\vec Q$, so both products are recomputed. For \texttt{XNysTrace}, the Lanczos iterations for $M\X$ are simply continued from the first run, so the refinement only costs the additional iterations. In our experiments we use $\tau_0 = 10^{-2}$ and $c = 10^{-2}$.

\subsubsection{Evaluation of $g_m$}
\label{subsubsec:evaluation-gm}
Each inner run of~\cref{subsubsec:inner-stopping} requires evaluating $g_m$ at the eigenvalues of its recurrence matrix, i.e., at the $r$ eigenvalues that define $\psi_r^{GL}[g_m,\vec P]$, see~\cref{eq:full-size-GL}, and at those of $T_j^{(i)}$ in~\cref{eq:inner-stopping}. By~\cref{thm:errorformula}, $g_m(z)$ is a divided difference of $f$ on the doubled Ritz values $\theta_1,\theta_1,\ldots,\theta_m,\theta_m$ and $z$. Evaluating this divided difference directly is prone to cancellation in floating-point arithmetic, especially when some of the nodes are close to each other. However, if $f$ is given by, or is closely related to, an integral of the form~\cref{eq:CauchyStieltjes}, then $g_m$ also admits the integral representation~\cref{eq:contourintegralerror}, which can be evaluated accurately by numerical quadrature. We now discuss such numerical quadrature rules for the functions used in our numerical experiments.

With $\varphi_m(t) \coloneqq (\beta_2\cdots\beta_{m+1})/\det(tI-T_m) = \prod_{j=1}^m \beta_{j+1}/(t-\theta_j)$, the representation~\cref{eq:contourintegralerror} reads
\begin{equation}\label{eq:gm-integral-phi}
	g_m(z) = \int_\Gamma \frac{g(t)\,\varphi_m(t)^2}{t-z}\,\mathrm{d}t.
\end{equation}
For all functions below, we evaluate a real parametrization of~\cref{eq:gm-integral-phi} by adaptive Gauss--Kronrod quadrature, as implemented in the Julia package QuadGK.jl~\cite{QuadGK}, for which the integration interval and breakpoints are specified for each function.

\paragraph{Stieltjes functions}
For a Stieltjes function~\cref{eq:stieltjes}, we have
\begin{equation}\label{eq:gm-stieltjes}
	g_m(z) = \int_0^\infty \frac{\rho(t)\,\varphi_m(-t)^2}{z+t}\,\mathrm{d}t.
\end{equation}
If $A$ is positive definite, then $\theta_j > 0$, and for $z \ge 0$ the integrand has the constant sign of $\rho$ and decays like $\rho(t)\,t^{-2m-1}$, so~\cref{eq:gm-stieltjes} can be evaluated without cancellation. We integrate over $[0,\infty)$ with breakpoints at $\theta_1,\ldots,\theta_m$. Examples are
\begin{equation*}
	x^{-\alpha} = \frac{\sin(\alpha\pi)}{\pi}\int_0^\infty \frac{t^{-\alpha}\,\mathrm{d}t}{x+t}, \quad 0<\alpha<1,
	\qquad\text{and}\qquad
	\frac{\log(1+x)}{x} = \int_1^\infty \frac{t^{-1}\,\mathrm{d}t}{x+t}.
\end{equation*}
Since divided differences of order $2m$ annihilate $\poly_{2m-1}$,~\cref{eq:gm-stieltjes} remains valid if $f(x) = \int_0^\infty \big[\rho(t)\,(x+t)^{-1} + q_t(x)\big]\,\mathrm{d}t$ with $q_t \in \poly_{2m-1}$, even if the two terms are not separately integrable. An important example is the entropy function~\cite{BenziRinelliSimunec2023}, for which
\begin{equation*}
	-x\log x = x(1-x)\int_0^\infty \frac{\mathrm{d}t}{(t+x)(t+1)}
	= \int_0^\infty \left[-\frac{t}{x+t} + \Big(1-\frac{x}{1+t}\Big)\right]\mathrm{d}t,
\end{equation*}
i.e., $\rho(t) = -t$ and $q_t \in \poly_1$. Hence, for every $m \ge 1$,
\begin{equation*}
	g_m(z) = -\int_0^\infty \frac{t\,\varphi_m(-t)^2}{z+t}\,\mathrm{d}t,
\end{equation*}
whose integrand again has constant sign. Similarly, $\log x = \int_0^\infty \big[(1+t)^{-1} - (x+t)^{-1}\big]\,\mathrm{d}t$ corresponds to $\rho(t) = -1$.

\paragraph{Exponential function $f(x) = e^{\beta x}$}
Here we use case~(i) with $g(s) = e^{\beta s}/(2\pi i)$ and $\Gamma$ the circle $s(t) = c + R e^{it}$, $t \in [0,2\pi]$, where $c = (a+b)/2$ and $[a,b]$ contains the Ritz values $\theta_j$ and all evaluation points $z$, e.g., $[a,b] = [\lambda_{\min},\lambda_{\max}]$. Since $\theta_j$ and $z$ are real, the integrand at $s(2\pi - t) = \overline{s(t)}$ is the complex conjugate of the integrand at $s(t)$, and~\cref{eq:gm-integral-phi} becomes
\begin{equation*}
	g_m(z) = \frac{1}{\pi}\int_0^\pi \operatorname{Re}\left[\frac{e^{\beta s(t)}\,\varphi_m(s(t))^2\,R e^{it}}{s(t)-z}\right]\mathrm{d}t.
\end{equation*}
The integrand is smooth on $[0,\pi]$, so no breakpoints are needed. Furthermore, we choose $R = (b-a)/2 + 2m/\abs{\beta}$, so that $\Gamma$ lies at distance $d = 2m/\abs{\beta}$ from $[a,b]$. This choice balances the growth of $e^{\beta s}$ against the decay of $\varphi_m(s)^2$ along $\Gamma$, since $d \mapsto e^{\abs{\beta} d} d^{-2m}$ is minimized at $d = 2m/\abs{\beta}$, and thereby limits cancellation in the quadrature. This is analogous to the choice of the optimal radius for computing high-order derivatives of analytic functions by Cauchy integrals, see Bornemann~\cite{Bornemann2011}; for $e^z$ and the $n$th derivative, this optimal radius is $n + \mathcal{O}(1)$.

\paragraph{Absolute value $f(x) = \abs{x}$}
This function is not analytic at $0$ and fits neither case~(i) nor case~(ii). However, for real $x$,
\begin{equation*}
	\abs{x} = \frac{2}{\pi}\int_0^\infty \frac{x^2}{x^2+u^2}\,\mathrm{d}u = \frac{2}{\pi}\int_0^\infty \Big[1 + u \operatorname{Im}\,(iu - x)^{-1}\Big]\,\mathrm{d}u,
\end{equation*}
i.e., up to a constant term, $\abs{x}$ is an integral of resolvents with poles $\pm iu$ on the imaginary axis. The constant is annihilated by the divided difference, and taking the imaginary part commutes with divided differences on real nodes, so
\begin{equation*}
	g_m(z) = \frac{2}{\pi}\int_0^\infty u\,\operatorname{Im}\left[\frac{\varphi_m(iu)^2}{iu-z}\right]\mathrm{d}u.
\end{equation*}
We integrate over $[0,\infty)$ with breakpoints at $\abs{\theta_1},\ldots,\abs{\theta_m}$ and $\abs{z}$. In contrast to the previous cases, the integrand does not have constant sign. Moreover, when $A$ is indefinite, $f$ is not smooth on $[\lambda_{\min},\lambda_{\max}]$, so the Gauss--Radau bounds of~\cref{sec:GaussRadau} and~\cref{prop:variance-error-operator} do not apply.

\section{Numerical experiments}
\label{sec:numerical-experiments}

In this section we compare our methods \texttt{GLHutch}, \texttt{GLXTrace} and \texttt{GLXNysTrace} against other trace estimators. 
Throughout the section $s$ denotes a budget parameter for the stochastic estimators and corrections: Hutchinson, \texttt{GLHutch}, and \texttt{GLXNysTrace} use $s$ test vectors, whereas \texttt{XTrace} and \texttt{GLXTrace} use $s/2$ test vectors. This corresponds to a budget of $s$ matrix-vector products with $f(A)$ or the error matrix $M$, excluding additional work due to the adaptive stopping criteria discussed in \cref{subsubsec:inner-stopping} and counting quadratic forms in Hutchinson and \texttt{GLHutch} in the same way as matrix-vector products. All numerical experiments are performed in Julia 1.12.6 on a laptop equipped with an AMD Ryzen 7 PRO 7840U processor and 16 GB of RAM.

\subsection{Log-determinant of the \texttt{bcsstk29} graph Laplacian}
\label{subsec:exp-bcsstk29-logdet}

We approximate $\log\det(A) = \trace(\log(A))$ for the shifted graph Laplacian $A = L + 10^{-1} I$, where $L = D - W$ is the Laplacian of the graph defined by the sparsity pattern of the matrix \texttt{HB/bcsstk29} from the SuiteSparse Matrix Collection~\cite{DavisHu11}, i.e., $W$ is its adjacency matrix and $D$ the diagonal matrix of the vertex degrees. Then $n = 13\,992$ and $\sigma(A) \subset [10^{-1}, 75.6]$. We compare the full-size global Lanczos approximation $\psi_m^{GL}$ for $m = 1,\dots,7$, Hutchinson's estimator with $s = 8,16,\dots,4096$ Rademacher vectors, \texttt{XTrace} with $s = 8,16,\dots,512$, and the hybrid methods \texttt{GLHutch}, \texttt{GLXTrace} and \texttt{GLXNysTrace} of~\cref{alg:gl-correction} with $m = 7$, for $s = 4,32,64,\dots,1024$ (\texttt{GLHutch}) and $s = 4,16,32,\dots,256$ (\texttt{GLXTrace} and \texttt{GLXNysTrace}). The quadratic forms of Hutchinson and \texttt{GLHutch} are approximated by global Lanczos, stopped by~\cref{eq:inner-stopping-hutch} with $C = 10^{-2}$; the block products of \texttt{XTrace}, \texttt{GLXTrace} and \texttt{GLXNysTrace} by block Lanczos, stopped by~\cref{eq:inner-stopping} with $c = \tau_0 = 10^{-2}$. The reference value $\log\det(A)$ is computed exactly from a sparse Cholesky factorization of $A$.

\begin{figure}[tbhp]
	\centering
	\includegraphics{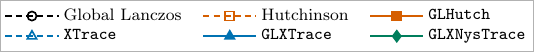}\par\medskip
	\subfloat[Normalized eigenvalue magnitudes]{%
		\includegraphics{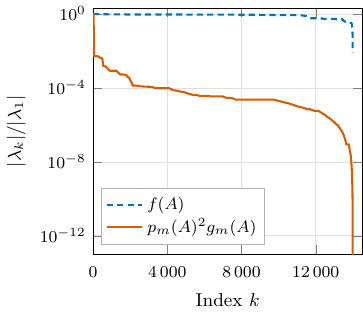}%
		\label{fig:timing-bcsstk29-logdet-spectra}}%
	\hfill
	\subfloat[Relative error against CPU time]{%
		\includegraphics{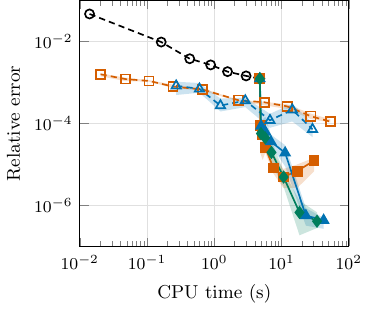}%
		\label{fig:timing-bcsstk29-logdet-timing}}
	\caption{Approximation of $\log\det(A)$, with $A$ the shifted graph Laplacian of \texttt{bcsstk29} ($n = 13\,992$). (a) Eigenvalues $\lambda_k$ of $\log(A)$ and of $p_m(A)^2g_m(A)$ for $m = 7$, sorted by decreasing absolute value and normalized by the largest one. (b) Relative error against CPU time. Markers show the median over $10$ runs and shaded bands indicate the central $50\%$ of the error distribution. Dashed lines with open markers denote global Lanczos and the stochastic estimators applied to $\log(A)$; solid lines with filled markers of the same colour and shape denote the corresponding hybrid method with $m = 7$ global Lanczos iterations, for budgets $s = 0,4,16,32,\dots,256$ ($s = 0,4,32,64,\dots,1024$ for \texttt{GLHutch}).}
	\label{fig:timing-bcsstk29-logdet}
\end{figure}

\Cref{fig:timing-bcsstk29-logdet} shows the relative error against the CPU time. Each stochastic estimator is run $10$ times for each $s$; we report the median relative error and the mean CPU time. The time of a hybrid method includes the $m = 7$ global Lanczos iterations, so its curve starts at the last point of the global Lanczos curve ($s = 0$).

Global Lanczos reduces the relative error from $4.7\cdot 10^{-2}$ to $1.2\cdot10^{-3}$ in $7$ iterations and $4.8$ seconds, but each iteration is more expensive than the previous one because of the fill-in of the basis matrices, and the convergence slows down. Hutchinson's estimator converges roughly at the expected rate $O(1/\sqrt{s})$, and reaches a relative error of $1.1\cdot 10^{-4}$ in $53$ seconds; \texttt{XTrace} is only slightly more accurate for the same time, with a relative error of $7.2\cdot10^{-5}$ in $29$ seconds. Already with $s = 4$, each of the three corrections reduces the error of $\psi_7^{GL}$ by more than an order of magnitude, to between $5.7\cdot10^{-5}$ and $8.8\cdot10^{-5}$, in at most $0.2$ additional seconds: after about $5$ seconds, the hybrid methods are already more accurate than Hutchinson after $53$ seconds, and comparable to \texttt{XTrace} after $29$ seconds. \texttt{GLHutch} stagnates between $5\cdot10^{-6}$ and $10^{-5}$ for $s \geq 128$, while the errors of \texttt{GLXTrace} and \texttt{GLXNysTrace} keep decreasing quickly with $s$: \texttt{XTrace} and \texttt{XNysTrace} are much more effective on the error matrix $M$ than \texttt{XTrace} is on $\log(A)$. This is consistent with the faster initial decay of the eigenvalue magnitudes of $M$ shown in \cref{fig:timing-bcsstk29-logdet-spectra}. \texttt{GLXNysTrace} reaches a relative error of $6.8\cdot 10^{-7}$ in $18$ seconds and $4.1\cdot10^{-7}$ in $34$ seconds, more than two orders of magnitude below the most accurate results of Hutchinson and \texttt{XTrace}; \texttt{GLXTrace} attains the same accuracy but requires about $25\%$ more time.

\subsection{Entropy of the Washington, DC road network Laplacian}
\label{subsec:exp-dc-entropy}

We approximate $\trace(f(L))$ with $f(x) = -x\log x$ (and $f(0) = 0$) for the graph Laplacian $L = D - W$ of the road network of Washington, DC from the 9th DIMACS Implementation Challenge~\cite{DIMACS9}. We take the graph unweighted and undirected, and restrict it to its largest connected component, which has $n = 9\,522$ vertices, $14\,807$ edges and maximum degree $6$; then $\sigma(L) \subset [0, 8.10]$. Since $\tau \coloneqq \trace(L) = 2 \cdot 14\,807$, the von Neumann entropy of the density matrix $\rho = L/\tau$ follows as $-\trace(\rho\log\rho) = \trace(f(L))/\tau + \log\tau$. The matrix $f(L)$ is indefinite, since $f > 0$ on $(0,1)$ and $f < 0$ on $(1,\infty)$, so that \texttt{XNysTrace} cannot be applied to it. The error matrix $M$, however, is negative semidefinite, since $g_m(z) < 0$ for $z \geq 0$ by the integral representation of~\cref{subsubsec:evaluation-gm}, so that \texttt{GLXNysTrace} can be applied to $-M$. We compare the full-size global Lanczos approximation $\psi_m^{GL}$ for $m = 1,\dots,8$, Hutchinson's estimator with $s = 8,16,\dots,4096$, \texttt{XTrace} with $s = 8,16,\dots,512$, and the hybrid methods with $m = 8$, for $s = 4,32,64,\dots,2048$ (\texttt{GLHutch}) and $s = 4,16,32,\dots,256$ (\texttt{GLXTrace} and \texttt{GLXNysTrace}). The Gauss--Radau node is $\lambda_{\min} = 0$, and the inner iterations are stopped as in~\cref{subsec:exp-bcsstk29-logdet}. The reference value $\trace(f(L)) = -4.0471\cdot 10^{4}$ is computed from the eigenvalues of $L$.

\begin{figure}[tbhp]
	\centering
	\includegraphics{timing-legend.pdf}\par\medskip
	\subfloat[Normalized eigenvalue magnitudes]{%
		\includegraphics{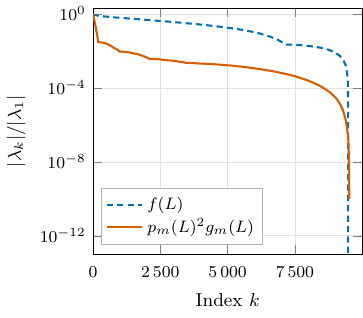}%
		\label{fig:timing-dc-entropy-spectra}}%
	\hfill
	\subfloat[Relative error against CPU time]{%
		\includegraphics{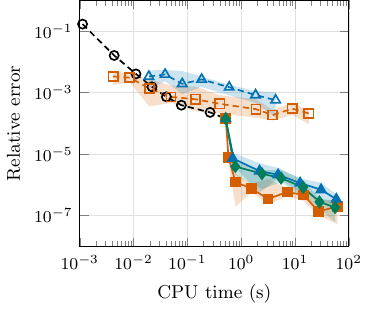}%
		\label{fig:timing-dc-entropy-timing}}
	\caption{Approximation of $\trace(-L\log L)$, with $L$ the graph Laplacian of the Washington, DC road network ($n = 9\,522$). (a) Eigenvalues $\lambda_k$ of $f(L)$ and of $p_m(L)^2g_m(L)$ for $m = 8$, sorted by decreasing absolute value and normalized by the largest one. (b) Relative error against CPU time. Markers show the median over $10$ runs and shaded bands indicate the central $50\%$ of the error distribution. Dashed lines with open markers denote global Lanczos and the stochastic estimators applied to $f(L)$; solid lines with filled markers of the same colour and shape denote the corresponding hybrid method with $m = 8$ global Lanczos iterations, for budgets $s = 0,4,16,32,\dots,256$ ($s = 0,4,32,64,\dots,2048$ for \texttt{GLHutch}).}
	\label{fig:timing-dc-entropy}
\end{figure}

\Cref{fig:timing-dc-entropy} shows the relative error against the CPU time, with the same statistics as in~\cref{fig:timing-bcsstk29-logdet}. Since every vertex has at most $6$ neighbours, the fill-in of the basis matrices grows slowly, and global Lanczos reduces the relative error from $1.7\cdot10^{-1}$ to $1.4\cdot10^{-4}$ in $8$ iterations and $0.52$ seconds. The stochastic estimators applied to $f(L)$ are much less effective: Hutchinson's estimator reaches a relative error of $1.9\cdot10^{-4}$ after $3.9$ seconds ($s = 1024$), but does not improve on this for $s = 2048$ and $4096$, and \texttt{XTrace} is less accurate than Hutchinson for the same time, with a relative error of $5.7\cdot10^{-4}$ after $4.4$ seconds ($s = 512$). The eigenvalues of $f(L)$ decay slowly in absolute value, with more than half of them within a factor $6$ of the largest, so that the low-rank part of \texttt{XTrace} does not pay off.

With $s = 4$, each of the three corrections reduces the error of $\psi_8^{GL}$ by a factor between $18$ and $36$, to between $4.0\cdot10^{-6}$ and $7.9\cdot10^{-6}$, in at most $0.26$ additional seconds; after $0.6$ seconds, \texttt{GLHutch} is about $25$ times more accurate than Hutchinson after $18$ seconds. In contrast to~\cref{subsec:exp-bcsstk29-logdet}, \texttt{GLHutch} is the most efficient hybrid method here: it reaches a relative error of $1.2\cdot10^{-6}$ in $0.78$ seconds ($s = 32$), $3.5\cdot10^{-7}$ in $3.1$ seconds ($s = 128$) and $1.4\cdot10^{-7}$ in $27$ seconds ($s = 1024$). For the same budget $s$, \texttt{GLXTrace} and \texttt{GLXNysTrace} are at most a few times more accurate than \texttt{GLHutch}, but their corrections are ten to twenty times more expensive, so that they need about $30$ seconds for a relative error of $3\cdot10^{-7}$ to $7\cdot10^{-7}$, and reach $3.4\cdot10^{-7}$ (\texttt{GLXTrace}) and $1.8\cdot10^{-7}$ (\texttt{GLXNysTrace}) only after about $57$ seconds ($s = 256$). Although the eigenvalues of $M$ decay faster in absolute value than those of $f(L)$, their initial decay is slower than in \cref{fig:timing-bcsstk29-logdet}, and the accuracy gains of the low-rank correction methods do not compensate for their higher computational cost.

\section*{Conclusions}
We developed a novel method for estimating the trace of a matrix function that combines a full-size global Lanczos approximation with a stochastic correction. As a byproduct of our analysis, we derived new a posteriori error bounds for global Lanczos approximations of any block size. Using tail bounds for stochastic estimators, we established convergence guarantees and illustrated the convergence behavior in numerical experiments. The numerical comparisons show that the new method can significantly outperform existing methods when they struggle to reach high accuracy. We also observed that the relative performance of the stochastic correction methods varies with the initial decay of the singular values of the global Lanczos error remainder. Explaining this empirical relationship and predicting which correction works best for a given problem are directions for future work.

\bibliographystyle{siamplain}
\bibliography{references}

\end{document}